\documentclass[11pt,a4paper]{article}

\usepackage[utf8]{inputenc}
\usepackage{amsmath}
\usepackage{amsfonts}
\usepackage{amssymb}
\usepackage{makeidx}
\usepackage{graphicx}
\graphicspath{ {./} }
\usepackage{listings}
\usepackage{tocloft}
\usepackage{verbatim}
\usepackage{xcolor}

\usepackage{calc}
\usepackage{CJKutf8}
\usepackage[hidelinks]{hyperref}

\usepackage{epigraph}
\newtheorem{le[mm]a}[theorem]{Le[mm]a}

\newcommand{\qed}{\nobreak \ifvmode \relax \else
      \ifdim\lastskip<1.5em \hskip-\lastskip
      \hskip1.5em plus0em minus0.5em \fi \nobreak
      \vrule height0.75em width0.5em depth0.25em\fi}
      
       \usepackage[polutonikogreek,english]{babel}
    
\usepackage[normalem]{ulem}

\usepackage{tikz, tkz-tab, pgfkeys}
\usetikzlibrary{calc,decorations.pathmorphing,shapes}

    \newcommand*{\tg}[1]{\textgreek{#1}}

  \newcounter{exer}[section]

    \makeindex

 \title{Euler's work on spherical geometry:  \\
 An overview with comments}
 \author{Athanase Papadopoulos and Vladimir Turaev}

\begin{document}

\maketitle

\maketitle

\vspace*{\fill} \epigraph{\itshape 
Wasser den Rhein hinunterflie\ss en, bis die Schule endlich merkt, da\ss \ die Mathematik eine Geisteswissenschaft sein kann und da\ss \ die Sch\"uler \textsc{Euler} so gut verstehen k\"onnen wie \textsc{Plato} oder \textsc{Goethe}. (Andreas Speiser, \emph{Introduction to Euler's} Introductio in analysin infinitorum. Tomus secundus. In:
Leonhardi Euleri opera omnia Vol. 9, Lipsiae 1945, \emph{p. {\sc xix}}.)
}

 \vfill\eject

\begin{abstract}

We review Euler's work on spherical geometry. After an introduction concerning the general place that trigonometric formulae occupy in geometry, we start by the two memoirs of Euler on spherical trigonometry, in which he establishes the trigonometric formulae using different methods, namely, the calculus of variations in the first memoir, and classical methods of solid geometry in the other.  In another memoir, Euler gives several formulae for the area of a spherical triangle in terms of its side lengths (these are ``spherical Heron formulae"). He uses this in the computation of numerical values of the solid angles of the five regular polyhedra, which is his goal in his memoir. We then review memoirs in which Euler systematically starts by establishing a theorem or a construction in Euclidean geometry and then proves an analogue in spherical geometry. We point out relations between Euler's memoirs on spherical trigonometry and works he did in astronomy, on the problem of drawing geographical maps, and  in geomagnetism. We also review some other works of Euler involving spheres, including a memoir on the three-dimensional Apollonius problem and others concerning algebraic curves on the sphere. Even though these works are not properly on spherical geometry, they show Euler's interests in various questions related to spheres and we think that they are worth highlighting in such an overview.  
Beyond spherical geometry,  the reader is invited to discover in this article
 an important facet of the work of the great Leonhard Euler. 
 This article will appear as a chapter in the book \emph{Spherical geometry in the eighteenth century, I:  Euler, Lagrange and Lambert}, Springer, 2026.

 \end{abstract}

\bigskip

 
  
  \bigskip
\tableofcontents

\section{Introduction}

We review Euler's memoirs on spherical geometry together with other works that include  his work on spherical geometry in a broad perspective. The plan of this article is the following.

We start (\S \ref{what-is}) with an excursus on the importance of trigonometric formulae, in any geometry, mentioning in particular the work of Lobachevsky and the paramount importance of the formulae of spherical trigonometry  that he established in his lectures and in several writings.

 After this digression which serves as one motivation	 for studying Euler's work on spherical trigonometry, we  review (\S \ref{s:trigo}) Euler's memoirs on this subject, namely, 
the
 \emph{Principes de la trigonom\'etrie sph\'erique tir\'es de la m\'ethode des plus grands et plus petits}
(Principles of spherical trigonometry deduced from the method of maxima and minima) \cite{Euler-Principes-T}, published in 1755, and  \emph{Trigonometria sphaerica universa, ex primis principiis breviter et dilucide derivata} (General spherical trigonometry, deduced from fundamental principles in a brief and clear manner)   \cite{Euler-Trigonometria-T}, published in 1782.
In each of these two memoirs, Euler establishes a set of  formulae of spherical trigonometry. 
 In the first memoir, he uses, in the proofs of these formulae, methods of differential geometry and the calculus of variations, whereas in the second he relies on more classical proofs that are based on the geometry of the 3-dimensional Euclidean space:  cutting by planes the solid angles defined by a spherical triangle and applying to the resulting Euclidean triangles the formulae of plane trigonometry. 
 
 In \S \ref{solid-angles}, we review Euler's 
  memoir \emph{De mensura angulorum solidorum} (On the measure of solid angles), in which he gives a geometrical proof of Girard's formula for the area of spherical triangles which is intrinsic to the sphere. The proof begins with a remark on the areas of lunes\index{lune (on the sphere)!area}  on the sphere\footnote{A \emph{lune} on the sphere is a figure bordered by two arcs of great circles.} that\index{lune (on the sphere)} goes back to Bonaventura Cavalieri. \index{Cavalieri, Bonaventura}  In the same memoir, Euler obtains several formulae for the area of a spherical triangle as a function of its side lengths, using analytical tools.  He then uses these formulae to obtain approximate values of the solid angles of the five regular polyhedra.  
  
  In \S \ref{s:Euler-locus}, we review Euler's memoir
  \emph{Variae speculationes super area triangulorum sphaericorum} (Various speculations on the area of spherical triangles) \cite{Euler-Variae-T}, in which he proves a result in spherical geometry which is the analogue of a well-known  result in Euclidean geometry, namely, a description of the locus of the vertices of triangles having a fixed base and a fixed area. In the Euclidean case, this locus is a straight line, and the result is contained in Euclid's \emph{Elements}. In the spherical case, the locus is not a geodesic, and the result is more complicated to prove. Euler attributes the problem and the result to his young collaborator Anders Johann Lexell, who wrote a memoir on this subject \cite{Lexell-Solutio}. 
  
Section \ref{Euclidean-to-spherical} is an overview of several memoirs in which Euler, like in the one which is the subject of the preceding section, proves some results in Euclidean geometry and then their analogues in spherical geometry. This includes the memoir \emph{Problematis cuiusdam Pappi Alexandrini constructio}
(A construction relative to a problem of Pappus of Alexandria) \cite{Euler-Pappi-T} 
 in which he obtains a generalization of a Euclidean construction from Pappus' \emph{Collection}, before considering the spherical case. This problem (in the Euclidean case) has attracted the attention of several mathematicians, in particular Euler's young collaborator at the Russian Academy,\index{Academy!Imperial Saint Petersburg Academy of Sciences} Nicolaus (better  known at the Saint Petersburg Academy as Nikolai Ivanovich) Fuss.\index{Fuss, Nikolai Ivanovitch}
The other mathematicians who worked on this problem include Lagrange, Castillon, Carnot, Brianchon and many others which we discuss in \S  \ref{Euclidean-to-spherical}. In another chapter of the present volume, Guillaume Th\'eret gives a detailed commentary on Pappus' problem.
In the same section, we review Euler's memoir \emph{Geometrica et sphaerica quaedam} (Some questions of the geometry of the plane and of the sphere) \cite{Euler-Geometrica-T}, in which he starts again by proving a result in Euclidean geometry and then its spherical analogue. This also gives us the occasion to review several works of Euler on Euclidean geometry.

In \S \ref{s:geography}, we discuss some of the spherical geometric arguments that appear in Euler's works on geography, astronomy and geomagnetism.

  In \S \ref{curves-on-the-sphere} and \ref{Apollonius4}, we survey some further work of Euler on the sphere (even though this work is not, strictly speaking, on spherical geometry), including his memoir on algebraic curves on the sphere (and other surfaces), a generalization to dimension three of Apollonius' problem of kissing circles (which becomes then  a problem of kissing spheres), making relations with recent problems, some of which are still part of current mathematical research.
  
In the memoirs we review in this article, beyond Euler's work on spherical geometry, it is Euler's work itself that the reader will discover.

We start with some thoughts on spherical trigonometry.

\section{Why spherical trigonometry?}\label{what-is}

The expression \emph{spherical trigonometry} refers to a set of formulae that make relations between the side lengths and the angles of a spherical triangle. Such formulae constitute a non-negligible part of Euler's mathematical heritage. In fact, even before writing his memoirs on this topic, Euler used extensively spherical trigonometry in his works on geography, astronomy and geomagnetism. But in our opinion, the importance of spherical trigonometry goes beyond these applications, as we shall explain later in this chapter, after including the question in a broader context.

First of all, let us state a general fact, namely, that the set of trigonometric formulae relative to a given geometry is an important aspect of that geometry, maybe the most important. Indeed, geometers know that the basic features of a geometry are contained in the properties of its triangles: the properties of the small (and infinitesimal) triangles tell us about the local properties of the geometry, and the properties of large triangles tell us about its large-scale properties.

 John Napier,\index{Napier, John} in his \emph{Mirifici logarithmorum canonis descriptio, ejusque usus, in utraque trigonometria; ut etiam in omni logistica mathematica, amplissimi, facillimi, \& expeditissimi explicatio} (1614)  \cite{Napier}, highlighted the importance of the trigonometric formulae. 
The following quote by Napier is borrowed from Annette A'Campo's chapter in the present volume \cite{AAC}: 
\begin{quote}\small
Since geometry is the art of measuring proposed quantities with complete accuracy, and as a figure [..] can be decomposed into the sides of triangles, or by further triangulation into smaller triangles, then a figure is composed of triangles with some of its angles and sides measured, and all the other triangular parts are to be found, and from which the extent of the figure can be found. Therefore it is clear that the arithmetical solution of any geometrical question depends on the principles by which triangles are solved.
\end{quote}

  It is interesting to remember in this respect that the great Lobachevsky,\index{Lobachevsky, Nikolai Ivanovich} when he started establishing the bases of the geometry that carries now his name, was very careful in establishing a complete set of trigonometric formulae. Less well known is that Lobachevsky, who was very much concerned about the question of the non-contradiction of his  geometry, argued that if there were a contradiction in his geometry, then there would be a contradiction stemming from the associated trigonometric formulae, for the reason we have just explained. But he had proved his (hyperbolic) trigonometric formulae using solely the geometry of the Euclidean plane and of the sphere. More precisely, he used the geometry of the horospheres of three-dimensional hyperbolic space (having shown beforehand that these horospheres naturally carry a Euclidean geometry) and that of the geometric spheres of his geometry (having shown beforehand that they bear the same geometry as those of the ordinary spheres of the three-dimensional Euclidean space). But then, because of the way he produced the trigonometric formulae of hyperbolic geometry, he was able to claim that if there were a contradiction stemming from these hyperbolic trigonometric formulae, then the contradiction would be already inherent in the formulae of Euclidean or spherical geometry.  Obviously, the question of the non-contradiction of spherical geometry was not addressed in those days, since this was not a geometry defined axiomatically, but a geometry that was part of Euclidean space geometry. Thus, Lobachevsky, without having a Euclidean model of hyperbolic geometry, reduced the non-contradiction problem of his geometry to that of Euclidean geometry. This fact is little known. Let us quote the remarks that he made on this question, at the end of his \emph{Elements of geometry} (the translation of the following text is from \cite[p. 223]{Rosenfeld1}):  
 
 \begin{quote}\small
  After we have found Equations (17)\footnote{Equations (16) and (17) of Lobachevsky's treatise to which he refers in the present passage are respectively sets of trigonometric formulae of spherical and hyperbolic geometry he established in his text, before this passage.} which represent the dependence of
the angles and sides of a triangle; when, finally, we have given the general
expressions for elements of lines, areas and volumes of solids, all else in
the Geometry is a matter of analytics, where calculations must necessarily
agree with each other, and we cannot discover anything new that is not
included in these first equations from which must be taken all relations of
geometric magnitudes, one to another.
  Thus if one now needs to assume
that some contradiction will force us subsequently to refute the principles
that we accepted in this geometry, then such contradiction can only hide in
the very Equations (17). We note, however, that these equations become
Equations (16) of spherical trigonometry as soon as, instead of the sides
 $a,b,c$ we put $a\sqrt{-1},b\sqrt{-1},c\sqrt{-1}$; but in ordinary Geometry and in
spherical Trigonometry there enter everywhere only ratios of lines; therefore
ordinary Geometry, Trigonometry and the new Geometry will always
agree among themselves.
  \end{quote}

It is also good to know that the trigonometric formulae (spherical and hyperbolic) are already
contained in Lobachevsky's first memoir, the \emph{Elements of geometry} \cite{Lobachevsky-Elements}, and in his
later works, such as the \emph{Geometrische Untersuchungen} \cite[\S 35 and 36 respectively]{Lobachevsky-Geometrische}.

Talking about the importance of the geometry of triangles, let us also recall that a large portion of the geometrical part of Euclid's \emph{Elements} and of  Menealus' \emph{Spherics} \cite{RP2} is a study of Euclidean, respectively spherical triangles (even if none of these treatises contains trigonometric formulae). 
Let us also recall that each of Herbert Busemann and A. D. Aleksandrov and M. Gromov developed  geometric theories that are based on a property of triangles, see \cite{Busemann-1948},  \cite{Aleksandrov-III} and \cite{Gromov}.
Finally, let us recall that W. Thurston, in his \emph{Geometry and topology of Three-manifolds} \cite{Thurston-GT3}, which arguably was the most influential 20th century work on geometry and topology, starts by establishing a set of trigonometric formulae for spherical and hyperbolic geometry.

   \section{On Euler's work on spherical trigonometry} \label{s:trigo}
  Euler wrote two memoirs on spherical trigonometry, \emph{Principes de la trigonom\'etrie sph\'erique tir\'es de la m\'ethode des plus grands et plus petits} (Principles of spherical trigonometry deduced from the method of maxima and minima) \cite{Euler-Principes-T}  and \emph{Trigonometria sphaerica universa, ex primis principiis breviter et dilucide derivata} (General spherical trigonometry, deduced from fundamental principles in a brief and clear manner) \cite{Euler-Trigonometria-T}. They are translated in Chapter 11 and 12 respectively of the present volume.
    Regarding these memoirs,  
 we start by quoting Rosenfeld, from his \emph{History of non-Euclidean geometry} \cite[p. 31]{Rosenfeld1}:
 
 \begin{quote}\small
 The modern form of spherical trigonometry, as well as of all trigonometry, is due to the great Leonhard Euler (1707-1783), a native of Basel, who worked in Petersburg and Berlin. Whereas trigonometry before Euler was concerned with trigonometric \emph{lines}, Euler's trigonometry dealt with trigonometric \emph{functions}, which he linked to the exponential function by means of the well known\index{sinus@\emph{sinus totus}} \emph{sinus totus},\footnote{In ancient sine tables, the \emph{sinus totus}  (the Latin expression for ``total sine") appeared as the radius of the circle that was used to construct this table. It corresponds there to the value of $\sin 90^{\mathrm{o}}$.} the complete sine, that is, the radius of the circle, and replaced it with unity.
 
 \end{quote}

To describe in some detail Euler's two memoirs, let us first recall that a spherical triangle\index{spherical!triangle} is a figure on the sphere consisting of three points (the vertices) together with three shortest segments (the sides) joining them pairwise.
Thus, a spherical triangle\index{triangle!spherical} gives rise to  six quantities: the three angles at the vertices and the three side lengths. 
By the usual criteria of equality of spherical triangles, the knowledge of any three of these six quantities determines the triangle. A trigonometric formula\index{trigonometric formula} is typically a relation between four or more of these quantities.  In each of his two memoirs, Euler derives a set of trigonometric formulae that can be useful in various situations, involving a wealth of variables like sums or differences of pairs of sides or of angles, etc.

   We first consider Euler's memoir \emph{Principes de la trigonom\'etrie sph\'erique tir\'es de la m\'ethode des plus grands et plus petits} \cite{Euler-Principes-T}, published in 1755.

 In the introduction to this memoir, Euler addresses the question of why one should use the \emph{methods of maxima and minima}, that is, the calculus of variations,\index{calculus of variations} to establish trigonometric formulae, while simpler methods exist. He answers this question in three points: 

(1) First, the methods of maxima and minima, he says, acquire a new domain of application in the field of spherical trigonometry. Indeed, it seems 
that this is the first time that trigonometric formulae are established using this method.\footnote{It should be recalled that the methods of the calculus of variations, since their introduction by Euler and then Lagrange, have been especially useful in the solution of problems of physics, even if these methods were formulated as part of mathematical analysis.
Euler's memoir \emph{Methodus inveniendi lineas curvas maximi minimive proprietate gaudentes, sive solutio problematis isoperimetrici latissimo sensu accepti} (A method for finding curved lines enjoying properties of maximum or minimum, or solution of isoperimetric problems in the broadest accepted sense), published in 1744 \cite{Euler-Methodus-T} and in which he gives a general presentation of the problems  that are dealt with using the calculus of variations as well as an exposition of the general methods for solving them, contains, by way of illustration, a list of one hundred special cases to which the general methods apply. These problems include a solution to the brachistochrone\index{brachistochrone problem} problem, others related to the principle of least action, others in geometric optics, and many more, all originating in physics. Talking about the calculus of variations, Euler writes in his memoir \emph{Principes de la trigonom\'etrie sph\'erique tir\'es de la m\'ethode des plus grands et des plus petits} \cite{Euler-Principes-T} that ``[\ldots]
 since we have shown that most of the mechanical and physical problems can be resolved very promptly with the aid of this method, it cannot be but very pleasant to see that the same method brings such a great help for the resolution of problems of pure Geometry."}

  (2) Euler declares that it is always useful to have different proofs of the same result, and that this will help for a better understanding of this result.

(3) Euler states that the methods of maxima and minima, unlike the previously known methods used for deriving the trigonometric formulae, are  very  general and can be applied to arbitrary surfaces.

Concerning the third point, one may recall that in the memoir  \emph{\'El\'ements de la trigonom\'etrie sph\'ero\"\i dique tir\'es de la m\'ethode des plus grands et des plus petits} (Elements of spheroidal trigonometry
drawn from the method of the maxima and minima) \cite{Euler-Elements-T},\footnote{The memoir is translated into English in the volume \cite{Caddeo-Papadopoulos} which is concerned by mathematical geography in the eighteenth century.} published the same year (1755) as the memoir \emph{Principes de la trigonom\'etrie sph\'erique tir\'es de la m\'ethode des plus grands et plus petits}, Euler studies triangles on a spheroid,\footnote{A spheroid\index{spheroid} is a surface obtained by the rotation of an ellipse along one of its axes. The study of the trigonometry of the spheroid was motivated by questions of geography. Indeed, at Euler's time, it was realized, based on Newton's work, that the Earth's shape was spheroidal rather than spherical. Therefore, spheroidal trigonometry became useful in cartography, an activity in which Euler was heavily involved. Euler writes in the same memoir
 \cite{Euler-Principes-T}: ``[\ldots] Since the surface of the Earth, being not spherical, but spheroidal, a triangle\index{triangle!spheroidal} formed on the surface of the Earth will belong to the species of which I have just talked about [that is, on the spheroid]."  The memoir  \cite{Euler-Elements-T} also contains explicit computations of distances between points on a spheroid.} and that he even starts there  by considering triangles on an arbitrary surface. Let us quote him from the introduction to that memoir (translation from \cite[p. 206]{Caddeo-Papadopoulos}):
\begin{quote}\small
Having established the Elements of spherical trigonometry on the principle of
the maxima and minima, my main goal was to fix such a general principle, from
which one could extract the resolution of triangles formed not only on a spherical
surface, but in general on an arbitrary surface. Since the sides of a spherical triangle
are arcs of great circles which, being the shortest paths that one can draw on the
surface of a sphere from one point to another, it is on an equal footing that I consider
the sides of a triangle described on an arbitrary surface, in such a way that they
become the shortest paths that lead from an angle to another on this surface. 
\end{quote}
 
Let us return to Euler's memoir \emph{Principes de la trigonom\'etrie sph\'erique etc.}  \cite{Euler-Principes-T}.  Euler starts   by recalling the definition of  a spherical triangle.\index{spherical!triangle} 
The idea of using the calculus of variations\index{calculus of variations}  in this setting is based on the fact that an edge of such a triangle is a geodesic on the sphere, that is, a shortest path between its endpoints. Finding the length of a side of a spherical  triangle becomes a variational problem. The first result in Euler's paper is the solution of this problem in the case where the side whose length is sought is opposite to a right angle. More precisely, he gives a solution to the following problem:

\bigskip
\noindent{\bf Problem.} Given the arc $AP$ on the equator  and the arc $PM$ on a meridian, to find on the spherical surface the shortest line $AM$ that can be drawn from point $A$ to point $M$.

\bigskip

The use of geographical terminology in this statement  is not unexpected, since Euler had in mind the example of the Earth, generally identified with a sphere in the problems of drawing geographical maps in which he had been heavily involved; see \cite{Caddeo-Papadopoulos}, and also \S \ref{s:geography} of the present chapter.

As in several  other memoirs of Euler, the main results in \cite{Euler-Principes-T} are presented in the form of problems and corollaries.   
  Euler starts by establishing a list of 30 formulae involving only three variables, valid for the case of a right triangle: Given any two quantities (sides or angles), he gives three formulae giving the values of the remaining ones (excluding the right angle, which is known).  Using this, he obtains a large collection of trigonometric  formulae for general (not necessarily right) triangles, involving sides,  angles and various combinations of them (half-sums of sides, differences of sides, etc.). Among them are formulae for the value of a side of a spherical triangle in terms of the three angles. Obviously, such formulae  cannot exist in Euclidean geometry, where the angles do not determine the sides.  
The memoir concludes with a proof using differential and integral calculus of the  formula attributed to Albert\index{Girard theorem}  Girard\footnote{\label{f:Girard1} Albert Girard (1595-1632) was a French-born mathematician who studied and spent all his life in the Netherlands, where his family was exiled for being protestant. He worked on algebra and number theory.  Like several mathematicians of his times, he was also an engineer. He wrote treatises on fortifications. He translated into French works of the  Flemish mathematician and music theorist Simon Stevin\index{Stevin, Simon} (1548-1620) and of the Dutch astronomer and mathematician Willebrord\index{Snellius, Willebrord} Snellius (1580-1626). His work is situated, historically, between those of Fran\c cois Vi\`ete\index{Vi\`ete, Fran\c cois} and Ren\'e Descartes.\index{Descartes, Ren\'e} Girard died at age 37.
Tannery wrote a biographical note on Girard \cite{Tannery}, in which he declares that the latter was one of the most remarkable geometers of the beginning of the 17th century.}
 which gives the area of a spherical triangle in terms of its angle excess.     
 
 Beyond the obtention of trigonometric formulae, Euler developed in the memoir \cite{Euler-Principes-T}, using the calculus of variations\index{calculus of variations}, a version of the intrinsic differential geometry of surfaces  embedded in 3-space, a theory that Gauss and, after him, Riemann, continued to develop in the century that followed.  Let us also mention in this respect that in the memoir
\emph{De linea brevissima in superficie quacunque duo quaelibet puncta iungente}
 (On the shortest line joining two points on a surface)
\cite{Euler-Brevissima-T} published in 1729, Euler (who was 22 years old) established the equation of a geodesic on an arbitrary surface.

In the memoir \emph{Trigonometria sphaerica universa, ex primis principiis breviter et dilucide derivata} (General spherical trigonometry, deduced from fundamental principles in a brief and clear manner), published in 1782
  \cite{Euler-Trigonometria-T}. Euler recovers the trigonometric formulae in a more classical manner, by methods that use the ambient three-dimensional space and with no differential or integral calculus.  He starts by establishing the following three basic formulae,
\[\frac{\sin C}{\sin c}=\frac{\sin A}{\sin a},
\]
\[
\cos A \sin c = \cos a \sin b - \sin a \cos b \cos C,
\]
and
\[
\cos c= \cos a \cos b +\sin a\sin b\cos C,
\]
from which he derives  the rest of the trigonometric formulae by performing algebraic transformations.
  Right at the beginning of the memoir, Euler notes that the choice of the notation $A,B,C$ and $a,b,c$ for the angles and sides opposite to them leads to a set of symmetric formulae. He highlights a rule for obtaining new formulae out of known ones, and this rule turns out to be a form of duality theory in spherical geometry, or the theory of the polar triangle.\index{polar triangle} This discovery was attributed by Michel Chasles,\index{Chasles, Michel} the famous mathematician and historian of mathematics, in \cite{Chasles-Apercu1}, to Fran\c cois\index{Vi\`ete, Fran\c cois} Vi\`ete (1540-1603), but it is known now that this discovery was already known to the Arab mathematicians of the 10th century; see the historical exposition in Chapter 10 of the present volume \cite{polar-history}.
       
 Several authors after Euler worked out, like him, a complete system of trigonometric formulae based on very few starting formulae. We mention in particular Lagrange, in his memoir \emph{Solution de quelques probl\`emes relatifs aux triangles sph\'eriques avec une analyse compl\`{e}te de ces triangles} (Solution of some problems relative to spherical triangles with a complete analysis of these triangles), published in 1799 and translated in the present volume  \cite{Lagrange-Solution-T},  and de Gua \cite{Gua}, whose work on the subject came shortly after that of Euler and before that of Lagrange.
 We end this section by quoting Delambre,\index{Delambre, Jean-Baptiste Joseph} from his \emph{Rapport historique sur les
progr\`es des sciences
math\'ematiques depuis 1789} (Historical report on the progress of mathematical sciences since 1789) \cite[p. 45]{Delambre-Rapport-H}, who includes Euler's work in a broader historical context:
 \begin{quote}\small
 Euler, in a memoir printed in 1779,\footnote{Delambre refers here to the memoir
\emph{Trigonometria sphaerica universa, ex primis principiis breviter et dilucide derivata}, published in 1782. The volume of the \emph{Acta} in which it appeared carries the mention ``pro anno 1779, pars prior", which explains Delambre's assertion on dates.} had\index{Delambre, Jean-Baptiste Joseph} reduced spherical trigonometry to a fully analytical form. Since everything is determined in a triangle of which we know, for example, two sides and the included angle, it follows that the whole of trigonometry is contained in the equation which serves to solve this particular case.  To extract it is a purely analytical problem. However, it requires some skill; it is not even useful to know in advance what one is looking for. This problem, moreover, was only a game for Euler. Mr. Lagrange, in treating the same subject in the sixth \emph{Cahier} of the \emph{Journal de l'\'Ecole Polytechnique}, added some curious propositions concerning the area of the spherical triangle and the pyramid for which this surface serves as a base. It should be noted, however, that this new way of considering trigonometry added almost nothing to the formulae we had in our possession, and we should not be surprised: the astronomers who made continual use of these methods had turned them upside down in every way; they had even arrived at solutions which are astonishing in their symmetry and elegance, when we consider the elementary means by which they had reached there. 
The research of Euler and Lagrange made the analytical method as elementary as the old one; and since that time several geometers, Messrs. Bertrand, Lacroix, Puissant and some others have varied this development, which in fact can be explained in several ways, starting with whichever of the particular formulae one wishes: given only one, the rest will follow from it.
 \end{quote}
 
 For the works of Bertrand, Lacroix and Puissant mentioned by 
 Delambre, see \cite{Bertrand, Lacroix, Puissant}.
 
To end this section, let us note that Euler, Lagrange and Lambert, in their memoirs on spherical trigonometry that are published in the present volume, all worked out sets of formulae that are suitable for computations with logarithms, that is, involving multiplications rather than additions. This was important for computations and for the establishment of tables that were used in astronomy and geography.
  
  \section{On the measure of solid angles}\label{solid-angles}

\begin{figure}
\centering
\includegraphics[width=0.8\linewidth]{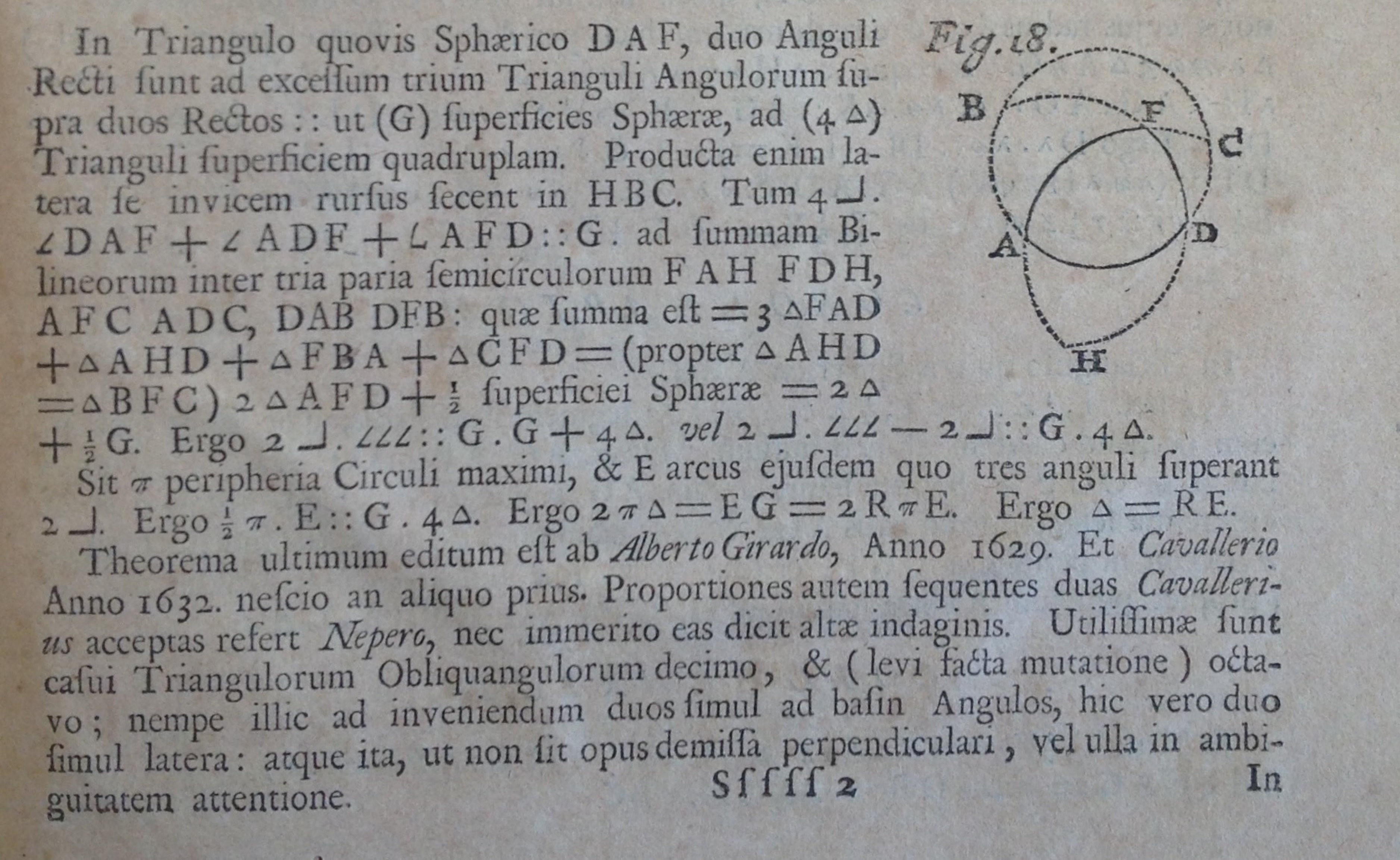}
\caption{From Wallis' proof of Girard's Theorem (\cite{Wallis2}, Vol.\ 2, p.\,875), the construction for the proof of Girard's formula. The region FAHDF is a lune.}   \label{fig:Wallis}
\end{figure}

\begin{figure}
\centering
 \includegraphics[width=0.8\linewidth]{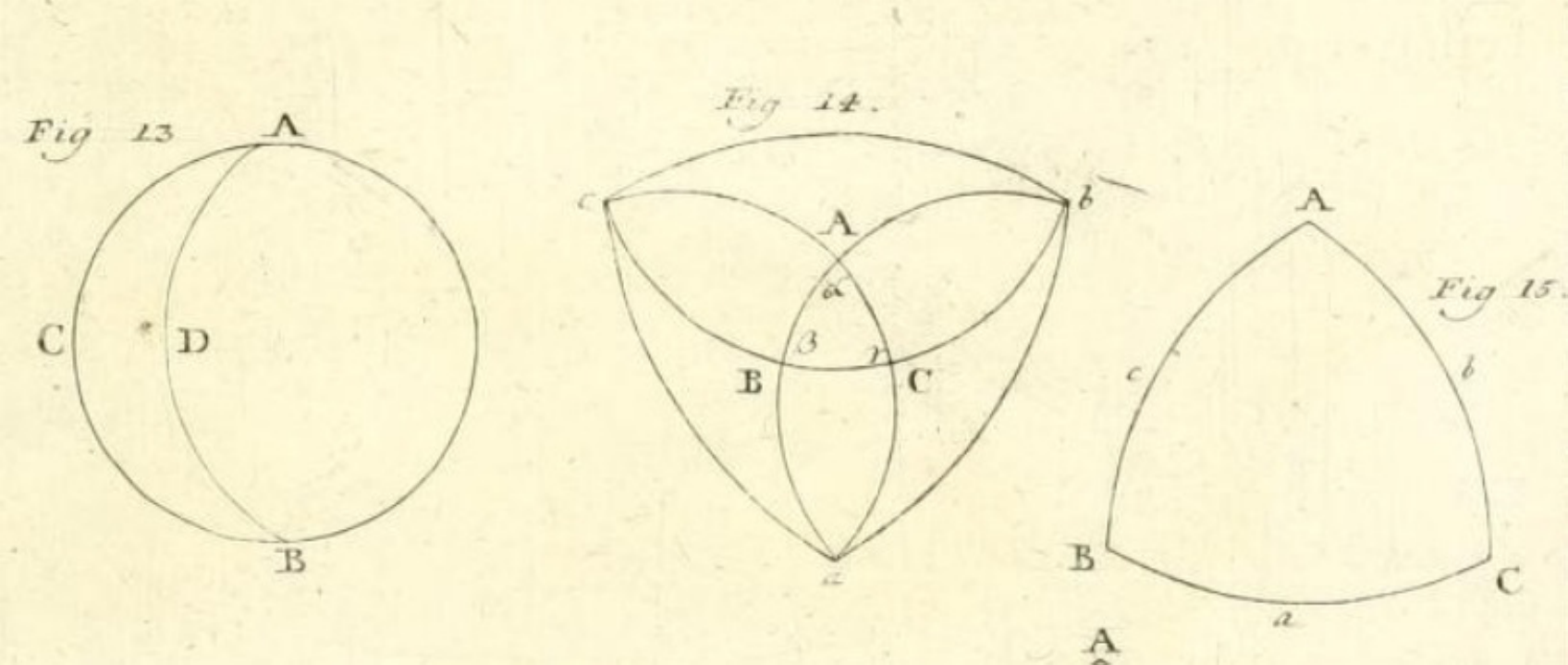}
\caption{From Euler's memoir \cite{Euler-Mensura-T}, the construction for the proof of Girard's formula. The region ACBDA (left) is a lune. In the middle figure, we have three lunes.}   \label{polar5}
\end{figure}

Spherical area is a central theme of Euler's memoir \emph{De mensura angulorum solidorum} (on the measure of solid angles) (1781) \cite{Euler-Mensura-T}. In this memoir, Euler's goal is to compute the numerical values of the solid angles of the five regular polyhedra. But 
the measure of a solid angle is the area of the region of the sphere of radius 1 centered at the vertex of this solid angle.\index{Girard theorem} This is the 3-dimensional analogue of the fact that the measure of an angle in the plane whose vertex is at the center of a circle of radius 1 is the length of the arc of circle subtended by this angle.
Thus, the computation of the solid angles becomes a question of computation of area in 
 spherical geometry. The simplest case is that of a solid angle formed by joining pairwise along their edges three Euclidean angles.  In this case, the value of the solid angle is equal to the area of a spherical triangle. Euler recalls in the memoir considered that, according to a\index{Girard, Albert} Girard's formula, the area of a spherical triangle is equal to its angle excess with respect to two right angles. He starts by giving a very simple proof of Girard's formula, based on the computation\index{lune (on the sphere)!area} of areas of lunes.\footnote{A \emph{lune},\index{lune (on the sphere)} in spherical geometry, is a region on the sphere bounded by two half-great circles, that is, a lune is bounded by two geodesic segments joining a pair of antipodal points, see Figure \ref{fig:Wallis} or Figure \ref{polar5}.} A similar proof was already given by John Wallis\index{Wallis, John} (1616-1703), who attributes it to Bonaventura Cavalieri\index{Cavalieri, Bonaventura} (1598-1647). See Figures \ref{fig:Wallis} and \ref{polar5}, extracted respectively from Wallis' and Euler's memoirs. The drawings summarize the main ideas of the proofs and their similarity is an indication of the fact that the two proofs are the same.\footnote{This proof with intersections of lunes can be adapted to hyperbolic geometry to prove the analogous result, namely, that the area of a hyperbolic triangle is equal to its angle deficit. Such a proof, in the hyperbolic case, was given by Gauss in \cite[Vol. VIII, p. 221]{Gauss-Works}.  Thurston, in his book \cite[p. 83 ff]{Thurston-book}, gives the same proof. The analogy between the Cavalieri--Wallis and the Gauss--Thurston proofs is accounted for in the article \cite{A-Campo-Papadopoulos}.} The excerpt from Wallis' proof, reproduced in Figure \ref{fig:Wallis}, mentions the names of Girard and Cavalieri.

Regarding the so-called Girard formula, let us note additionally that Lagrange, in his memoir \emph{Solution de quelques probl\`emes relatifs aux triangles sph\'eriques avec une analyse
compl\`ete de ces triangles}, which is translated in the present volume, Chapter 18 \cite{Lagrange-Solution-T}, writes the following: 
``[\ldots] we should rather attribute this theorem to Cavalieri\index{Cavalieri, Bonaventura} who gave it in the \emph{Directorium generale uranometricum}, printed in Bologna in 1632, with the beautiful proof reported on by Wallis and inserted since then in most of the Trigonometries."

 After the proof he gives in \cite{Euler-Mensura-T} of the so-called Girard formula, Euler establishes several formulae for the area of a spherical triangle in terms of its side lengths. This is a spherical analogue of the Euclidean Heron formula. 
This allows Euler to give numerical values for  the solid angles of all five regular polyhedra, after he solves the following problem:   

\emph{Given a solid angle consisting of $n$ planar angles which are equal to a certain quantity $a$ and which are equally inclined among themselves, to find the measure of the solid angle.}

It is conceivable that this is the first time in history where numerical values of the solid angles of the five regular polyhedra are computed.
 Incidentally, after Euler gives the spherical Heron formula, he shows that it  reduces to the usual (Euclidean) Heron formula when the three sides of the spherical triangle become infinitely small.

\section{The locus of spherical triangles with a given area} \label{s:Euler-locus}

In the memoir \emph{Variae speculationes super area triangulorum sphaericorum} (Various speculations on the area of spherical triangles) \cite{Euler-Variae-T}, published in 1797, Euler starts by giving several formulae for the area of a spherical triangle\index{triangle!spherical} as a function of its side lengths. He then considers the problem of finding the locus of the vertices of a spherical triangle with a given base and fixed area. The analogous problem, in the Euclidean case, is the subject of Propositions 37 and 39 of Book I of Euclid's \emph{Elements}. In this case, the desired locus is a straight line. At the same time, Euclid proves that any triangle in this family of triangles can be cut into pieces that can be joined together to form a square having the same area as the initial triangle. We recall that the \emph{Elements} do not contain any definition of area, but only the notion of ``two figures having the same area". The fact that two figures can be cut into pieces that can be glued again to form the same figure says that the two initial figures have the same area.

Unlike the Euclidean case, in spherical geometry, the sought locus is not a straight line (i.e. a great circle), but a small circle, that is, a locus of points equidistant from a great circle. Furthermore, contrary to what one might expect by analogy with the Euclidean case, the great circle in question is not the one containing the base of the triangle.

 Euler writes in \S 16 of his memoir that the ``occasion for these reflections" was given to him by the theorem made known by his young collaborator\index{Academy!Imperial Saint Petersburg Academy of Sciences} Anders\index{Lexell, Anders Johan} Johan Lexell.\footnote{Anders Johan Lexell (1740-1784) started working at the Saint Petersburg Academy of Sciences under Euler's guidance  and he eventually became an important astronomer and geometer at the Academy. The reader may refer to the articles \cite{Zhukova}  by Alena Zhukova  and \cite{Papa-Inde2} by the first author for an account on Lexell's work on spherical geometry.}
 Lexell wrote a memoir on this subject, titled \emph{Solutio problematis geometrici ex doctrina sphericorum} (Solution of a geometric problem in spherical geometry) \cite{Lexell-Solutio}. This memoir is analyzed and commented on in the sequel of the present volume, \cite{Caddeo-AP-iii}.
 Sorlin and Gergonne,\index{Gergonne, Joseph Diez} in their 1825  article
\cite{Sorlin-Gergonne}, give a solution to what is called now Lexell's problem,\index{Lexell problem} without quoting Euler, but quoting Lexell and Legendre.

In 1855, Lebesgue published a paper  titled \emph{D\'emonstration du th\'eor\`eme de Lexell} (A proof of Lexell's theorem) \cite{Lebesgue1855} in which he gives a proof of that theorem which in fact is Euler's proof. 
 At the end of this paper, the editor of the journal adds a comment  saying that Eug\`ene Catalan,\index{Catalan, Eug\`ene} in his  \emph{\'El\'ements de G\'eom\'etrie} (Book VII, Problem VII), gives a proof of Lexell's theorem. Catalan does not mention Euler. 
 The only original contribution of Lebesgue in this paper consists in two remarks that he makes after the proof. The first one says that there is  an analogous result for surfaces of revolution which admit an equator. The second remark says that by taking the stereographic projection of the figure on the sphere, we get  a family of curved triangles in the Euclidean plane that have a given base, whose angle sum is constant, and whose vertices are on the same circle. 
 

\section{Other memoirs: from Euclidean to spherical geometry}\label{Euclidean-to-spherical}

We have already noted that the so-called Lexell problem,\index{Lexell problem} which Euler studies in spherical geometry, is the spherical version of the  result contained in Euclid's \emph{Elements} concerning the locus of the vertices of the triangles having the same base and a fixed area.  
 We also said that in his memoir on the measure of solid angles \cite{Euler-Mensura-T}, Euler gives a spherical analogue of the Euclidean
 Heron formula\index{Heron formula!spherical} for the area of a  triangle in terms of its side lengths.    There are two other memoirs by Euler, translated in the present volume, in which Euler obtains results in spherical geometry that are analogues of results in Euclidean geometry. We shall now discuss them.
  
In the memoir \emph{Problematis cuiusdam Pappi Alexandrini constructio}
(A construction relative to a problem of Pappus of Alexandria), published in 1783
 \cite{Euler-Pappi-T}, Euler solves,\index{Pappus problem} in the Euclidean and then in the spherical setting,\index{Pappus of Alexandria} a problem inspired by a construction of Pappus,\footnote{Pappus of Alexandria was a prominent mathematician who lived around the end of the third century and the first half of the fourth century. His major work, a substantial part of which has survived, is the \emph{Collection} (Synagoge, \tg{Sunagwg'h}), which contains mathematical results, some of which are difficult to prove, covering much of geometry.
His name is attached to a theorem concerning the alignment of three points (and which is therefore not the one Euler is talking about), which is at the foundations of projective geometry. Pappus also introduced the notion of anharmonic ratio. His \emph{Collection} contains a host of theorems with detailed proofs (sometimes several proofs), some of which are difficult, which are not due to him, but which would not have come down to us without him. Some are taken from works by Euclid, Archimedes, Heron, Apollonius and others which are completely lost.
The \emph{Collection} contains results in Euclidean geometry, spherical geometry, analytic geometry, astronomy, and mechanics, as well as problems of the type that would today be called ``recreational mathematics." Jean-\'Etienne Montucla
 in \cite[t. 1, p. 329]{Montucla} writes\index{Montucla, Jean-\'Etienne} that Pappus ``gives in his \emph{Mathematical Collection} marks of an uncommon intelligence
in geometry, and in several places one finds
traces of genius." Among the mathematicians who have immersed themselves in the study of Pappus, at least Ren\'e Descartes,.\index{Descartes, Ren\'e} Thomas Simpson\index{Simpson, Thomas} and Michel Chasles\index{Chasles, Michel} should be mentioned.
The latter published a work, which he called \emph{Les trois livres des porismes d'Euclide, r\'etablis pour la premi\`ere fois d'apr\`es la notice et les lemmes de Pappus et conform\'ement au sentiment de R. Simson\index{Simson Robert} sur la forme des \'enonc\'es de ces proposition} (The three books of Euclid's porisms recovered for the first time according to the notice and lemmas of Pappus and in accordance with the feeling of R. Simson on the form of the statements of these propositions) \cite{Chasles-Porismes}. The work is entirely drawn from Pappus's \emph{Collection}.
Pappus's \emph{Collection} is both a great source of inspiration for mathematicians and an invaluable document on the history of mathematics. The Introduction to Book VII of that work contains a significant amount of unique information on ancient mathematicians and their results.}
from the \emph{Collection} \cite[Proposition 117 of Book VII]{Pappus-Eecke}. The problem is to construct a triangle which is inscribed in a given circle and whose three sides are contained in three lines that are required to pass through three given points. The construction is non-trivial.  Pappus gives this construction in the Euclidean setting and only in the special case where the three given points are collinear.\footnote{Pappus' \emph{Collection} contains sections on spherical geometry, but the problem which Euler solves in both the Euclidean and the spherical case was considered by Pappus only in the Euclidean case.} 

This problem has attracted the attention of several prominent mathematicians. 
In the same volume of the  \emph{Acta} in which Euler's memoir was published, another solution of the same problem, in the spherical setting as well, was given by\index{Academy!Imperial Saint Petersburg Academy of Sciences} Nicolaus\index{Fuss, Nikolai Ivanovitch} Fuss,\footnote{Nicolaus (or Nikolai Ivanovich) Fuss (1755-1826), originating from Basel, was sent to Saint Petersburg by  Daniel Bernoulli to assist Euler in his works. He later became one of the most brilliant geometers there, with important contributions to spherical geometry. Fuss became later Euler's dear friend and he married his granddaughter, Albertine Benedikte Philippine Luise Euler, the daughter of his son Johann Albrecht, who was also a mathematician.} another\index{Euler, Johann Albrecht} young\index{Euler,  Albertine Benedikte Philippine Luise} collaborator of Euler, in a memoir titled \emph{Solutio problematis geometrici Pappi Alexandrini} (Solution of a geometrical problem of Pappus of Alexandria) \cite{Fuss-Solutio}.\footnote{An English translation of this memoir by Fuss together with commentaries is published in the forthcoming volume \cite{Caddeo-AP-iii}.} An elementary solution by Lagrange of the same problem, based on trigonometric formulae, appeared in the memoir  \emph{Solution alg\'ebrique d'un probl\`eme de g\'eom\'etrie}, whose English translation is contained in Chapter 17 of the present volume, see \cite{Lagrange0}.
 Michel Chasles,\index{Chasles, Michel} in his \emph{ Aper\c cu historique sur l'origine et le d\'eveloppement des m\'ethodes en g\'eom\'etrie} \cite[p. 328-329]{Chasles-Apercu1} has a note on the history of this so-called Pappus problem (in the Euclidean setting). According to this account, Euler's proof (in the Euclidean setting) can be traced back to a theorem attributed to the Scottish mathematician Matthew\index{Stewart, Matthew} Stewart concerning four points taken arbitrarily on a straight line \cite{Stewart}.     We also learn from the same source that in 1742, the Swiss mathematician Gabriel Cramer\index{Cramer, Gabriel} (1704-1752) proposed this problem in the general (Euclidean) setting of non-collinear points\index{Cramer--Castillon problem (or Castillon problem)} to Johann Castillon\index{Castillon, Johann (Giovanni Francesco Salvemini da Castiglione)} (1704-1791) and that the latter published a solution, in 1776, in the \emph{Nouveaux m\'emoires de l'Acad\'emie
royale des sciences et belles-lettres de Berlin}. It is interesting to note that both Cramer and Castillon were correspondents of Euler, the first one between 1743 and 1751, and the second one from 1745 to 1765. This correspondence is published in \cite{Euler-Corresp-Cramer-Castillon}. 
Chasles also notes that 
Annibale Giordano (1769-1835),\index{Giordano, Annibale Giuseppe Nicol\`o} a 16 year old Napolitan,  and\index{Malfatti, Gian Francesco}  Gian Francesco Malfatti (1731-1807)\footnote{Malfatti's name is attached to the so-called \emph{Malfatti problem}\index{Malfatti problem} which asks for the construction of three circles inside a triangle that are tangent to each other and to the sides of the triangle.} solved the more general problem in which a triangle is replaced by a polygon with an arbitrary number of sides passing through an equal number of points. Their works are published in Volume IV of the \emph{Memorie della Societa Italiana delle scienze detta dei XL}.  In the \emph{Berlin Memoirs} of 1796, the Swiss mathematician Simon l'Huillier\index{L'Huillier, Simon Antoine Jean}  (1750-1840) published a modified algebraic solution of the same general problem. Lazare\index{Carnot, Lazare} Carnot (1753-1823),  in his \emph{G\'eom\'etrie de position} \cite[p. 383]{Carnot}, gave a modified form of Lagrange's solution and applied it  to the case of an arbitrary polygon, a case already considered by Giordano and Malfatti. Chasles also mentions solutions to generalizations of Pappus' problem by  Charles-Julien\index{Brianchon, Charles-Julien} Brianchon (1783-1864) in the case where a conic section replaces the circle, and where the three points are aligned, by Joseph Diez\index{Gergonne, Joseph Diez} Gergonne  (1771-1859) in the general case considered by Brianchon and allowing a solution only with straightedge (and no compass), and by Jean-Victor\index{Poncelet, Jean-Victor} Poncelet (1788-1867), who, in his \emph{Trait\'e des propri\'et\'es projectives} \cite{Poncelet} considers an extension to the case of a polygon with an arbitrary number of sides and where the circle is replaced by an arbitrary conic. A detailed mathematical commentary on Euler's result 
 presented in this memoir, with various approaches and generalizations by other mathematicians, is provided by Guillaume Th\'eret in Chapter 4 of the present volume \cite{Theret-Castillon}.

  Now we pass to the memoir \emph{Geometrica et sphaerica quaedam} (Some questions of the geometry of the plane and of the sphere) \cite{Euler-Geometrica-T}, which also contains results in Euclidean plane geometry followed by their analogues in spherical geometry. In this memoir, Euler starts by proving the following theorem: 
    \medskip
 
 \emph{Let $ABC$ be a triangle in the plane and let $a,b,c$ be points on the sides $BC, AC,AB$ respectively.
If the lines $Aa, Bb,Cc$ intersect at a common point $O$, then we have
$$
\frac{AO}{Oa}\cdot\frac{BO}{Ob}\cdot\frac{CO}{Oc}=\frac{AO}{Oa}+\frac{BO}{Ob}+ \frac{CO}{Oc}+2.
 $$}
  
 He also gives the following converse:
   
 \medskip
\emph{In the Euclidean plane, given three segments $AOa, BOb, COc$ meeting at a common point $O$ and satisfying the above equation, we can construct a triangle $ABC$ such that the points $a,b,c$ are on the sides $BC, AC,AB$ respectively.}

In the Euclidean case, Euler derives several different proofs of the theorem and its converse. He also considers the case where the point $O$ lies outside the triangle. In the spherical case, he obtains an analogous result. The new condition is:\footnote{ To avoid dealing with special cases, we are assuming assume that the spherical triangle is contained in a quarter sphere.}
$$
 \frac{\tan AO}{\tan Oa}\cdot\frac{\tan BO}{\tan Ob}\cdot\frac{\tan CO}{\tan Oc}=\frac{\tan AO}{\tan Oa}+\frac{\tan BO}{\tan Ob}+ \frac{\tan CO}{\tan Oc}+2.
$$
Thus, to pass from the Euclidean to the spherical case, the equation satisfied is the same up to replacing the lengths by their tangents. Mathematicians are fond of such analogies, and they seek to understand the reasons for them. Euler gives two proofs of this spherical result, one of them  based on spherical trigonometry, and the second one using a projection of the triangle onto a plane which is tangent to the sphere at the point $O$, reducing the problem to the Euclidean one, which he considered before.

It seems that the questions studied in this memoir have not been studied before Euler, neither from the Euclidean nor from the spherical point of view.
In the paper \cite{2015} written with Weixu Su, the first author of the present chapter has worked out the hyperbolic case, where a similar condition holds, with the tangent function replaced by the hyperbolic tangent, that is, the condition reads:
$$
 \frac{\tanh AO}{\tanh Oa}\cdot\frac{\tanh BO}{\tanh Ob}\cdot\frac{\tanh CO}{\tanh Oc}=\frac{\tanh AO}{\tanh Oa}+\frac{\tanh BO}{\tanh Ob}+ \frac{\tanh CO}{\tanh Oc}+2.
$$

To end this section, and since we are talking about Euclidean geometry, let us mention that Euler wrote several memoirs  in which he obtained a large number of results in this domain. Among these memoirs, we mention the 
\emph{Solutio problematis geometrici circa lunulas a circulis formatas} Solution of a geometric problem on the lunules  made by circles) \cite{E73}, in which he gives a solution  to a problem formulated by Daniel Bernoulli\index{Bernoulli, Daniel} in his \emph{Exercitationes quaedam mathematicae} (Some mathematical exercises)  \cite{Bernoulli-Exercitationes} (1724),\footnote{Euler says that this problem was proposed by Christian Goldbach.}  related to the question of dissecting circles into pieces that have equal area. It is interesting to  read Euler's comments on the solution he presents: 
\begin{quote}\small [The solution] is not only presented geometrically without analysis, but analysis is also considered less suitable for its solution.
This disadvantage of analysis, although commonly resorted to in many geometric problems, nevertheless seems to me to be attributable less to analysis than to the analyst.
In this problem, I clearly prove that not only is analysis unsuitable for its solution, but that the geometric method is far preferable. Using this method, I will present a general solution to this problem, which can be obtained by a geometric method.
\end{quote}

One may compare this statement with the fact that a few years later, and concerning other kinds of problems, Euler would seek to replace geometric arguments with analytical ones, and it was in the person of Lagrange that he found someone who could do it. 
The latter's ideas simplified Euler's approach in that they avoided tricky geometric arguments by replacing them with analytical ones, which led directly to what later came to be called the Euler--Lagrange equation. Euler expressed this in a letter dated September 6, 1755, in which he wrote to Lagrange:
\begin{quote}\small
After reading your last letter, in which you seem to have brought the theory of \emph{maxima} and \emph{minima} almost to its highest point of perfection, I cannot sufficiently express my admiration for the extraordinary sagacity of your mind. For not only had I called for in my treatise on this subject a purely analytical method which would allow the rules set out therein to be deduced, but I had subsequently devoted much effort to the search for such a method: and you have truly given me great joy by having the kindness of sharing with me your very profound and very solid reflections on these questions; this is why I acknowledge my debt to you.
\end{quote}

Among Euler's other memoirs on Euclidean geometry, we mention the \emph{Solutio facilis problematum quorundam geometricorum difficillimorum} (Easy solutions to some difficult geometrical problems) \cite{E325} and
\emph{Variae demostrationes geometriae}
(Various geometric demonstrations)
\cite{E135}. In the latter, he proves several elementary results involving the perimeter, area and radius of inscribed circles and other characteristics of Euclidean figures. For example, one result says that the area of a Euclidean rectangle is equal to the product of half the sum of its side lengths with the radius of the inscribed circle.

In the memoir 
\emph{Solutio problematis geometrici}
 (Solution of a problem of geometry) \cite{E192}, Euler presents his solution to another  problem of Pappus,
namely, given two conjugate diameters of an ellipse, to find the conjugate axes. (This is Problem XVII of
Book VIII of the Collection, p. 851 of ver Eecke's edition \cite{Pappus-Eecke}.) The problem is also quoted  in T. L. Heath's, \emph{History of Greek Mathematics} \cite[Vol. II, p. 437-438]{Heath-History}.  
 Talking about Heath, let us also point out a note in his edition of  Euclid's \emph{Elements} \cite[Vol. 1 p. 402]{Heath} in which he quotes the following theorem of Euler:  \emph{In any quadrilateral the sum of the squares on the sides is equal to the sum of the
squares on the diagonals and four times the square on the line joining the middle points of the diagonals.}
  Euler published several other elementary (but not necessarily easy to prove) results on Euclidean geometry and it would be interesting to work out the analogous statements on the sphere and in the hyperbolic plane.

 \section{Geography, astronomy and geomagnetism}\label{s:geography}

  Spherical geometry has applications in geography, astronomy and other fields. Chasles writes in his \emph{Aper\c cu historique}
\cite[p. 235]{Chasles-Apercu1} that since the works of the ancients (he mentions Theodosius,
Menelaus and Ptolemy), ``if the theory [of spherical geometry] was extended and if it
has attained, in the hands of our most celebrated geometers, a high degree
of perfection, this was always done by almost preserving the same framework,
because the goal was always the same: the computation of triangles in order
to fit the service of the astronomer and the navigator, and for the great
geodesic operations which showed the true form of the Earth." Euler
was interested in all domains of applications of mathematics, in particular in geography, astronomy and geomagnetism, which all involve spherical trigonometry.

As already noted in \S \ref{s:trigo},  Euler, in his memoir \emph{Principes de la trigonom\'etrie sph\'erique tir\'es de la m\'ethode des plus grands et plus petits} \cite{Euler-Principes-T},  used  geographical  terminology while he established his trigonometric formulae. In the introduction of this memoir, he talks about using the same methods for the study of the spheroid, ``the surface of the Earth, being not spherical, but spheroidal", and announcing his forthcoming memoir 
\emph{\'El\'emens de la trigonom\'etrie sph\'ero\"\i dique tir\'es de la m\'ethode des plus grands et plus petits} \cite{Euler-Elements-T} in which he establishes, by the same methods,  the trigonometric formulae  for the surface of the spheroid and from which  he deduces explicit distances between locations on the Earth.

     The main question in the problem of drawing geographical maps is to find maps of least distortion (with an appropriate meaning for this word) from a subset of the sphere into the plane. The volume \cite{Caddeo-Papadopoulos} contains memoirs by Euler, Lagrange and Lambert on this theory,  translated into English and with commentaries.  In Euler's memoirs on this subject, the geometry of the sphere, and in particular spherical trigonometry, plays a major role. Astronomy is also linked to spherical geometry, because the celestial vault is assimilated to a sphere, and making a map of the sky and studying the movements of the stars amount to questions of spherical geometry.
In any case, Euler's memoirs on astronomy involve some amount of spherical trigonometry. 
For example,  his memoir \emph{Solutio problematis astronomici ex datis tribus stellae fixae altitudinibus et temporum differentiis invenire elevationem poli et declinationem stellae} (Solution to problems of astronomy: given the altitudes and time differences for three fixed stars, to find the elevation of the pole and the declination of the star) \cite{E14}, published long before his memoirs on spherical trigonometry, consists almost entirely of considerations of spherical trigonometry, with some numerical examples in the last part. The goal of this memoir is to solve the problem announced in the title. Euler
  starts by recalling the following result which he attributes to\index{Mayer, Friedrich Christoph}  Friedrich Christoph Mayer:\footnote{ Friedrich Christoph Mayer (1762-1841) was a mathematician who, like Euler,  was working at the  Imperial Academy of Sciences in Saint Petersburg. From the
time of Euler's arrival to this Academy, in 1727, Mayer helped him with
various problems of celestial mechanics applied to the determination of the sun's orbit,
the movement of the planets and the calculation of lunar eclipses. In addition, Mayer worked
with Jakob Hermann and Daniel Bernoulli on a theory
of the moon\index{theory of the moon} based on Delisle's extensive programme of astronomical observation, which
included the occultations of stars and planets by the moon, as well as lunar and solar
eclipses. At that time, Delisle was the major astronomer and geographer of the Academy.}
  In any spherical triangle $ABC$, we have 
\[\cos A=\frac{\cos BC-\cos AB \cdot \cos AC}{\sin AB\cdot \sin AC}
.\] 
He then proves the following statement, quoted as a theorem. 
  
   \emph{In any spherical triangle $ABC$, we have}

\[
\begin{split}
 \cos BC & =
  \frac{\cos(AB+AC)-\cos (AB-AC)}{2}\\
 & +\frac{\cos A\cdot \cos(AB-AC)-\cos A\cdot \cos (AB+AC)}{2}.
\end{split}
\]  
 
Likewise, the memoir 
  \emph{De inventione longitudinis locorum ex observata
lunae distantia a quadam stella fixa cognita} 
(On finding the longitude of a place by observing the distance between the moon and a known fixed star) \cite{E570} contains a substantial part on spherical trigonometry. There are many other examples.

On geomagnetism, we  mention Euler's memoir \emph{Recherches sur la d\'eclinaison de l'aiguille aimant\'ee} (Researches on the declination of the magnetized needle) \cite{E237}. The (magnetic) declination\index{magnetic declination} at a given point on the Earth's surface is the angle between the magnetic north and true north. In this memoir, Euler mentions the work of Edmund Halley, who proposed a new idea for showing the declination on geographical maps, see 
\cite{Halley}. The memoir \cite{E237} starts by the following: ``Since this research, as well as the following ones, requires the analytical resolution of spherical triangles, it will be good to put the formulae before our eyes, so that we do not need to look for them elsewhere. I will therefore begin with right triangles; etc." The memoir, which is 77 pages long, contains a complete set of formulae of spherical trigonometry (without proofs), that are used all along.

    Finally, let us mention that spherical geometry is the basis of an argument that is the subject of a letter sent by Euler to Maupertuis\index{Maupertuis, Pierre Louis Moreau de} dated July 4, 1744 (letter No. 7, p. 49 of \cite{Euler-Opera-IV-VI}). The argument is analyzed in Chapter 4, titled \emph{On the Duration of the Passage of a Star
from an Almucantar to Another} in volume \cite{Caddeo-Papadopoulos}, see \cite{Charitos-Mucantarat}.

 \section{Curves on the sphere, the spheroid, etc.}\label{curves-on-the-sphere}
 
 The works mentioned in this section are not on spherical geometry, but we have collected them here because they involve the sphere. Indeed, Euler studies in these memoirs curves on a sphere embedded in the 3-dimensional Euclidean space.
 
In the memoir \emph{De curva rectificabili in superficie sphaerica} (On rectifiable curves   on the surface of the sphere) \cite{Euler-Curva1771} (1771), Euler studies  curves on the sphere that are expressible by algebraic equations.\footnote{Euler calls these curves ``rectifiable", a word used today in a different meaning.} He says that this problem  arises from what he calls the ``great Florentine problem [\ldots] raised on a spherical surface a century ago." 
He finds such a curve, namely, the
``spherical epicycloid", obtained by the motion on a great circle of a small
circle whose diameter makes a rational ratio with the diameter of the sphere. 

From this question regarding the surface of the sphere, he moves on to more general surfaces.

In the  memoir  \emph{De lineis rectificabilibus in superficie sphaeroidica quacunque geometrice ducendis} (On the rectifiable lines that have to be drawn geometrically on a spheroidal surface) \cite{Euler-Lineis} published the same year, Euler gives a simpler proof of the result of the preceding one \cite{Euler-Curva1771}, at the same time generalizing it from the case of the sphere to the case of a spheroid. In the introduction, he points out the difficulty of   
the question of drawing  geometrically algebraic lines on a surface. He writes: ``A most difficult question is when it is required to draw geometrically lines on the surface of any body, with all of whose arcs defined algebraically; since in cylinders, besides straight lines parallel to the axis, no other rectifiable curves can be drawn."
 
In the memoir \emph{De curvis rectificabilibus in superficie coni recti ducendis} (On rectifiable curves drawn on the surface of a right cone)  \cite{Euler-coni}, published in 1781, Euler finds a family of rectifiable curves on a right cone whose ratio of side-length to the diameter of the base is rational. Let us quote from the introduction: \begin{quote}\small
Since we can certainly affirm that no rectifiable line can be drawn on a cylindrical surface except the straight lines themselves parallel to the axis of the cylinder, the same seems to be true of conical surfaces, because a cylinder can be viewed as a sort of cone, while the height is increased infinitely.
However, I have observed that there are such right cones on the surface of which, in addition to the straight lines drawn from the vertex of the cone to the circumference of the base, innumerable other curved lines can be geometrically described, admitting rectifications; but this does not seem to be possible on all right cones, much less on oblique ones.
\end{quote}

\section{Another problem on spheres}\label{Apollonius4}

\begin{figure}
\centering
\includegraphics[width=0.5\linewidth]{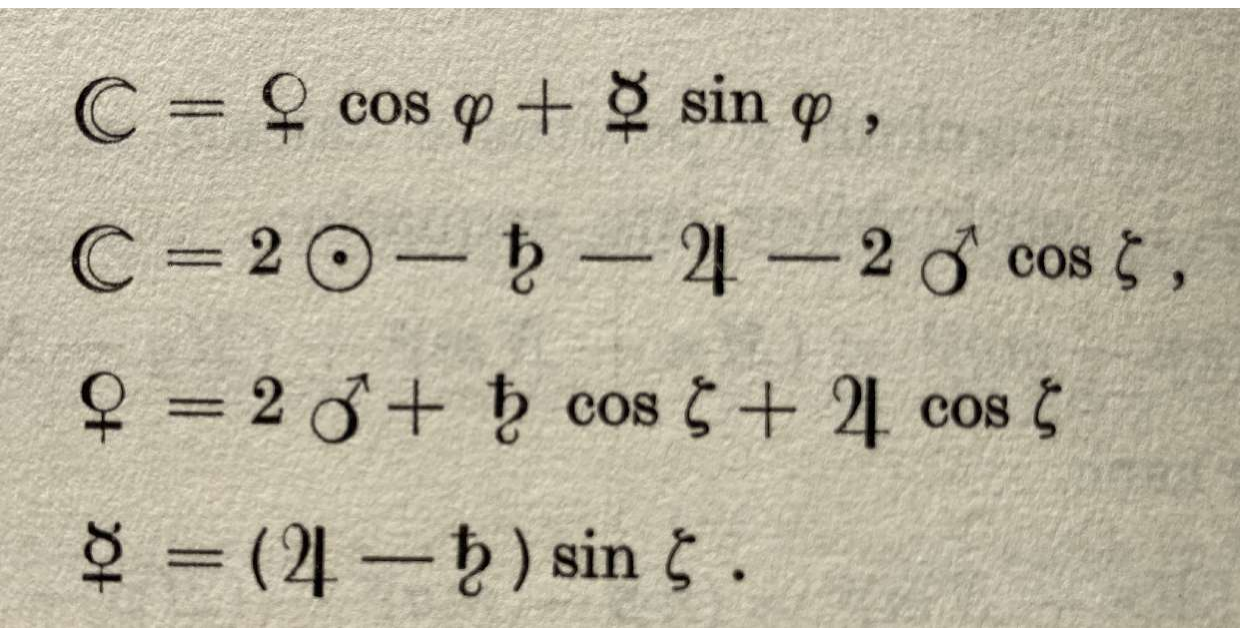}
\caption{A plate from Euler's memoir \emph{Solutio facilis problematis, quo quaeritur sphaera, quae datas quatuor sphaeras utcunque dispositas contingat} }  \label{fig:Zodiaque-Euler}
\end{figure}

In this section, we  consider Euler's memoir \emph{Solutio facilis problematis, quo quaeritur sphaera, quae datas quatuor sphaeras utcunque dispositas contingat} (An easy solution to the problem of finding a sphere that touches four given spheres, arranged in an arbitrary way) \cite{E733}. This is a three-dimensional analogue of the famous Apollonius' problem in the plane.\footnote{The reference is to Apollonius of Perga\index{Apollonius of Perga} (ca. 262-190 BC).  The problem, in its simplest form, is to construct (using straightedge and compass) a circle that is tangent to three given circles in a plane.
Apollonius' problem is part of his  work  \tg{'epafa'i} (``Contacts", or ``Tangencies"), which does not survive but which is reported on by Pappus in his \emph{Collection}. A detailed content of what is known about this work of Apollonius, and a history of the attempts made for its reconstruction, are given by Ver Eecke, the Belgian editor of Pappus'  \emph{Collection}, see \cite[Introduction, p. LXVI ff.]{Pappus-Eecke} We learn from there that Apollonius' work on this problem presumably contained twenty-one lemmas, sixty theorems and eleven problems, which Pappus, in his \emph{Collection}, summarizes in a unique proposition saying that given three arbitrary elements which could be points, lines or circles, to construct a circle which passes by the points (in the case where the objects are circles) 
or is tangent to the lines or circles. Pappus enumerates ten cases in which Apollonius solved the problem.} Like the problems considered in the preceding section, strictly speaking, Euler's memoir is not on spherical geometry, nevertheless it concerns spheres. 
 Euler presents two solutions to the problem of constructing a sphere that is tangent to four given spheres. The proofs rely on trigonometric formulae.
Towards the end of this memoir, Euler begins using astrological signs as variables, see Figure \ref{fig:Zodiaque-Euler}.\footnote{Lambert used these signs in his memoir on taxiometry,  see Chapter 6 of the present volume, and the figure reproduced there.}

At the end of the memoir, Euler writes:
%
%
\begin{quote}\small
In this way, we find a pair of solutions to a very difficult problem, which at first sight, seemed to require the most abstruse stereometric investigations, problems of this kind usually requiring both the most intricate figures and the most troublesome calculations, while the solutions given here can be obtained with the help of not too lengthy calculations.
The problem itself is not new, but was once found to have been solved by the great geometer Fermat;\index{Fermat, Pierre} but since at that time the calculation of angles was almost completely unknown, it is not surprising that our solution is found to be much more convenient.
\end{quote}

There is a relationship between the two-dimensional Apollonius' problem, which is at the origin of the problem discussed in this memoir, and the problem of Pappus, which is considered by Euler in the memoir \cite{Euler-Pappi-T} and on which we commented in \S \ref{Euclidean-to-spherical}. Guillaume Th\'eret mentions this relationship in his commentary on Euler's memoir, \cite{Theret-Castillon}, Chapter 4 of this volume. The reconstruction of Apollonius' proof and more generally the work on this problem attracted the attention of the greatest geometers such as Vi\`ete, Descartes, Newton, 
 Poncelet,  Gergonne, Cauchy and, as we said, Euler in his memoir \emph{Solutio facilis problematis, quo quaeritur circulus, qui datos tres circulos tangat} \cite{E648}.
 The three-dimensional version, treated by Euler in this memoir, was also considered by Fermat\index{Fermat@de Fermat, Pierre} in his \emph{De contactibus sphaericis}, see \cite[p. 52-69]{Fermat-Oeuvres}. 
 
 On this higher-dimensional Apollonius problem, besides Euler's memoir \cite{E648}, we mention the works of 
   Carnot \cite{Carnot}, Hachette \cite{Hachette1808, Hachette1815}, Binet \cite{Binet1815}, Battaglini \cite{Battaglini1851}, Fran\c cais \cite{Francais1810, Francais1813}, 
  Dupin \cite{Dupin1813}, Serret \cite{Serret1848}, Coaklay \cite{Coaklay1859-1860} and Alvord
  \cite{Alvord1882}.   Chasles, in his \emph{Aper\c cu historique sur l'origine et le d\'eveloppement des m\'ethodes en g\'eom\'etrie} \cite[p. 66]{Chasles-Apercu1}, writes about Euler's work on the three-dimensional Appollonius problem: 
\begin{quote}\small
The question of the sphere tangent to four others is one in which Geometry has long had the upper hand over Analysis. In 1779, Euler presented two analytical solutions of it to the Academy of Petersburg, which did not appear until the beginning of this
century, in the Volume of Academy for the years 1807-1808 (printed in 1810).
Carnot had already indicated an analytical solution in his \emph{G\'eom\'etrie de position} (p. 416),
but without carrying out the developments that would have led to a second-degree equation.
Today, it is M. Poisson who was the first to completely resolve this question by
calculation (Bulletin de la Soci\'et\'e philomatique, ann. 1812, p. 141.) Shortly afterwards,
MM. Binet and Fran\c cais also gave different analytical solutions (17th issue of the 
Journal de l'\'ecole polytechnique, and  tome III  of the Annales de Math\'ematiques).
\end{quote}

Naturally, Apollonius' problem in the plane has been extended to a similar problem on the sphere and on other quadratic surfaces. The question on the sphere is to construct all spherical circles that are tangent to three given circles, see \cite{Gergonne-Recherche, Carnot, Vansson}.

To conclude this section, let us say that Euler has studied the classical (2-dimensional) Apollonius problem in his memoir
\emph{Solutio facilis problematis, quo quaeritur circulus, qui datos tres circulos tangat} (An easy solution of a problem, in which a circle is searched for,  tangent to three given circles) \cite{E648},  written in 1779 and published in 1810.
This problem has been generalized to the case of $n>3$ circles and to the case  where  some of the circles are replaced by  lines or by conics, to the case of intersecting circles, etc. There is a very large bibliography on generalized Apollonius problems. Among the mathematicians who worked on it, we mention again  Pappus, then Descartes, Newton, Gauss, Cauchy, Poncelet, Gergonne, Laguerre, and there are many others.

Jean-Alfred Serret,\index{Serret, Jean-Alfred} in his paper \emph{De la sph\`ere tangente \`a quatre sph\`eres donn\'ees} (On the sphere tangent to four given spheres) \cite{Serret1848}, solves the problem of a sphere tangent to four given spheres using the notion of power of a point with respect to a sphere (the square of the distance from the center minus the square of the radius). Gaston Darboux\index{Darboux, Gaston}, in his paper \emph{Sur les relations entre les groupes de points, de cercles et de
sph\`eres dans le plan et dans l'espace} (On the relationships between groups of points, circles and spheres in the plane and in space) \cite{Darboux-relations} solves the problem of constructing a circle intersecting three given circles at given angles. In his paper, Darboux says that the problem has been proposed by Jakob Steiner.\index{Steiner, Jakob}

Among the modern works around generalized Apollonius' problems, we cannot help but think of the theorem of the Koebe--Andreev--Thurston theorem,\index{Koebe--Andreev--Thurston theorem} see \cite[Chapter 13]{Thurston-GT3} and the review in \cite{Bowers}.

Finally, talking about Euler's work on spherical geometry, we should mention the works on this topic done by the young mathematicians of the Saint Petersburg Academy of Sciences\index{Academy!Imperial Saint Petersburg Academy of Sciences} who worked under his direction or who were inspired by him. We have devoted a forthcoming volume to this subject, \cite{Caddeo-AP-iii}. The first author has also published a short overview of this topic \cite{Papa-Inde2}.

\noindent {\bf Authors' addresses}: 

\medskip

\noindent Athanase Papadopoulos, Institut de Recherche Mathématique Avancée (Université de Strasbourg et CNRS), 7 rue René Descartes, 67084 Strasbourg Cedex France,

 \noindent email : athanase.papadopoulos@math.unistra.fr 
 
 \medskip

\noindent Vladimir Turaev, Department of Mathematics, Indiana University, Rawles Hall, Bloomington, IN. 47405,  USA
 
 \noindent email: vtouraev@indiana.edu

\printindex

 \end{document}